\documentclass[11pt,reqno]{amsart}
\usepackage{amscd} 
\usepackage{amsfonts} 
\usepackage{amssymb} 
\usepackage{latexsym} 
\input{diagrams}
\newcommand{\ncm}{\newcommand}

\newtheorem{theorem}{Theorem}[section]
\newtheorem{prop}[theorem]{Proposition}
\newtheorem{lemma}[theorem]{Lemma}
\newtheorem{cor}[theorem]{Corollary}
\newtheorem{lem&def}[theorem]{Lemma \& Definition}

\newtheorem{example}[theorem]{Example}

\def\C{\mathbb{C}\,} 
\def\Z{\mathbb{Z}\,}

\def\|{\, |\,}

\ncm{\End}{\mbox{\rm End}\,}

\def\Hom{\mbox{\rm Hom}\,}
\def\Im{\mbox{\rm Im}\,}

\def\id{\mbox{\rm id}}

\def\into{\hookrightarrow}
\def\to{\rightarrow}
\def\Ind{\mbox{\rm Ind}}
\def\Res{\mbox{\rm Res}}

\def\o{\otimes}    
\def\bra{\langle}
\def\ket{\rangle}

\ncm{\rarr}[1]{\stackrel{#1}{\longrightarrow}}
\ncm{\larr}[1]{\stackrel{#1}{\longleftarrow}}
\def\cop{\Delta}

\def\eps{\varepsilon}

\def\du1{\hat 1}

\def\-1{_{(-1)}}
\def\0{_{(0)}}
\def\1{_{(1)}}
\def\2{_{(2)}}
\def\3{_{(3)}}

\def\du1{\hat 1}
\def\lact{\triangleright}
\def\ract{\triangleleft}

\begin{document}

\title[Depth Two, Normality and a Trace Ideal Condition]{Depth Two, Normality and a Trace Ideal 
Condition for Frobenius Extensions}
\author{Lars Kadison} 
\address{University of New Hampshire \\ Kingsbury Hall \\ 
Durham, NH 03824 U.S.A.} 
\email{kadison@math.unh.edu}

\thanks{}

\author{Burkhard K\"ulshammer}
\address{Mathematisches Institut\\
Friedrich-Schiller-Universit\"at \\
 07740 Jena \\
GERMANY}   
\email{kuelshammer@uni-jena.de}
\date{} 

\begin{abstract}
A ring extension $A \| B$ is  depth two
if its tensor-square satisfies a projectivity condition
w.r.t.\ the bimodules ${}_AA_B$ and ${}_BA_A$. In this case
the structures $(A \o_B A)^B$ and $\End {}_BA_B$ are bialgebroids
over the centralizer $C_A(B)$ and there is a certain Galois
theory associated to the extension and its endomorphism ring.
  We specialize the notion of depth two 
  to
induced representations of semisimple algebras and character theory of finite groups.  
We show that depth two subgroups over
the complex numbers are normal subgroups. As a converse, we observe that normal
Hopf subalgebras over a field are depth two extensions. 
A generalized Miyashita-Ulbrich action on the centralizer of a ring extension is introduced, 
and applied to a  study of depth two and separable extensions, which yields new characterizations of separable
and H-separable extensions.  With a view to the problem of when separable
extensions are Frobenius, we supply a trace ideal condition for when a ring
extension is Frobenius.  
\end{abstract} 
\maketitle

\section{Introduction and Preliminaries}

In noncommutative Galois theory we have the classical notions of Frobenius
extension and separable extension. For example, finitely generated Hopf-Galois extensions
are Frobenius extensions via a nondegenerate trace map 
defined by the action of the integral element $t$
in the Hopf algebra.  Also a Hopf-Galois extension is separable if its Hopf algebra
is semisimple or if the counit of $t$ is nonzero. 

There have been various efforts
to define a notion of noncommutative normal extension by Elliger and others. 
Recently a notion of depth
two for Bratteli diagrams of pairs of $C^*$-algebras has been widened to Frobenius
extensions in \cite{KN} and to ring extensions in \cite{KS} for the purpose
of reconstructing bialgebroids of various types depending on the hypotheses that
are placed on the ring extension and its centralizer.  
For example, a depth two
balanced Frobenius extension of algebras with trivial centralizer is a Hopf-Galois
 extension (with normal basis property), since dual Hopf algebras with natural actions
 are constructible on the step two centralizers
in the first levels of the Jones tower.    At the other extreme, a depth two
balanced extension has dual bialgebroids over the centralizer with Galois
actions.  
In between these two types of extensions we have extensions that
have natural Galois actions coming from Hopf algebroids or weak Hopf algebras
(e.g., groupoid algebras):
dual antipodes appearing when we place the Frobenius condition on a depth two extension.
 In Section~3, we extend the analogy of normal subfields and their Galois correspondence with normal subgroups
 to show that depth two subgroups of finite groups, 
over the complex numbers, are normal subgroups. We also observe a converse more generally stated
for normal Hopf subalgebras (that they are depth two extensions).  

There are intriguing problems in how to draw a Venn diagram for the various ring extensions
in noncommutative Galois theory.  For example, are split, separable extensions with conditions
of finite generation and projectivity automatically Frobenius extensions or quasi-Frobenius (QF) extensions
\cite{CK,KS}? Or an H-separable extension is separable and of depth two, but what
extra condition do we need for a separable depth two extension to be H-separable?  In Sections~2, 4, and 5
we work towards a clarification of these two questions by
focusing on a generalized Miyashita-Ulbrich action on the centralizer denoted by the module $R_T$.  
In Hopf-Galois theory (or normal subgroup theory) for a ring extension (or group ring extension)
$A \| B$ with
centralizer $R$, this is a right action defined 
by $$r \ract h := t^1 r t^2 := r \cdot t$$ where an element $h$ of the Hopf algebra $H$ corresponds to an
element $t := \beta^{-1}(1 \o h)$ in  $T : = (A \o_B A)^B$ under the Galois isomorphism $\beta: A \o_B A \stackrel{\cong}{\rightarrow}
A \o H$. The ring $T$ is also an $R$-bialgebroid in the depth two theory \cite{KS}.  
We will see in sections~3 and~4 that various conditions on the module $R_T$ thus defined and a
ternary product isomorphism involving this module lead to depth two, separable or
H-separable extensions. In section~5 we begin with a brief introduction to Frobenius extensions
and then create a small theory of trace ideal condition and Morita
context from the basic idea that an extension $A \| B$ should be Frobenius if there
is a bimodule map $E : A \to B$ and Casimir element $e = \sum_i e^1_i \o e^2_i$
such that $\sum_i E(e^1_i) e^2_i = 1$.

In this paper, a ring extension $A \| B$ is any unital ring homomorphism $B \to A$, called
proper if the map is monic, subrings being  the most important
case.  This induces a $B$-$B$-bimodule structure
on $A$ via
pullback denoted by ${}_BA_B$. Given this object of study we fix
a series of notations important to  the study of depth two.
  Denote its endomorphism
ring $S := \End {}_BA_B$.  Denote the group of homomorphisms
$\Hom ({}_BA_B, {}_BB_B) := \hat{A}$, a right  $S$-module with respect 
to ordinary composition of functions, 
or, if $A \| B$ is proper, identifiable with a right ideal in 
$S$. 

Recall that the centralizer $R$ of the ring extension $A \| B$ above is defined
as the $B$-central elements of ${}_BA_B$
 with notations $R := C_A(B) := A^B$. Then there are two ring mappings
of $R$ into $S$ given by $\lambda: R \to S$, $\lambda(r)(a) = ra$,
a ring homomorphism, and $\rho: R \to S$, $\rho(r)(a) = ar$,
a ring antihomomorphism in our convention.  Note that the two mapping
commute at all pairs of image points in $S$:
\begin{equation}
\label{eq: commute}
\lambda(r) \rho(r') = \rho(r') \lambda(r) \ \ \ \ (r,r' \in R)
\end{equation}
The two commuting maps induce two  $R$-$R$-bimodule structures on $S$:
a left $R^e$-structure 
denoted by ${}_{\lambda,\, \rho}S$ and given by ($\alpha \in S, r,r' \in R$)
\begin{equation}
\label{eq: left R^e structure on S}
{}_{\lambda,\, \rho}S:\ \ \ 
r \cdot \alpha \cdot r' =
 \lambda(r) \rho(r') \alpha =  r \alpha(-) r'
\end{equation}
and a right $R^e$-structure on $S$ denoted by $S_{\lambda,\, \rho}$
and given by 
\begin{equation}
\label{eq: right R^e structure on S}
S_{\lambda,\, \rho}:\ \ \  r \cdot \alpha \cdot r' := \alpha \lambda(r') \rho(r) = \alpha(r'(-)r).
\end{equation} 
We note that $S$ itself forms an $R^e$-$R^e$-bimodule under these two
actions.  

In the tensor-square $A \o_B A$ we draw the reader's attention to
the group of 
so-called $B$-casimir elements $T : = (A \o_B A)^B$ within.
$T$ has a natural ring structure induced from the isomorphism
$T \cong \End {}_AA\! \o_B\! A_A$ given by 
$$ t \longmapsto (a \o a' \mapsto ata')$$ for each $t \in T, a,a' \in A$
with inverse $F \mapsto F(1 \o 1)$:  the induced ring structure has unity
element $1_T = 1 \o 1$ with multiplication
\begin{equation}
\label{eq: T-multiplication}
tt' = {t'}^1 t^1 \o t^2 {t'}^2
\end{equation}
where we suppress a possible summation $t = \sum_i t^1_i \o t^2_i$.  
There is a natural $T$-$A^e$-bimodule on $A \o_B A$ stemming
from ``inner'' multiplication from the left  (as above) and ``outer'' multiplication
from the right.  We note the right $T$-module structure on 
the centralizer $R$ given by
$$ R_T: \ \ r \cdot t = t^1 r t^2 \ \ (r \in R, t \in T),$$
a cyclic module since $1 \o r \in T$ for every $r \in R$. 
Finally the two mappings $\sigma, \tau: R \to T$ given by $\tau(r) = r \o 1$,
an antihomorphism of rings, and $\sigma(r) = 1 \o r$, a ring homomorphism,
commute as in Eq.~(\ref{eq: commute}), 
whence they induce two commuting $R^e$-bimodule structures denoted by 
${}_{\sigma,\tau}T$ and $T_{\sigma,\tau}$ given
by 
\begin{equation}
\label{eq: outer structure on T}
T_{\sigma,\tau}: \ \ \ r \cdot t \cdot r' := t \sigma(r')\tau(r) = rt^1 \o t^2 r'
\end{equation}
which is the same as the restriction to $R$ of the ordinary bimodule ${}_AA \! \o_B \! A_A$,
and the more exotic bimodule
\begin{equation}
\label{inner structure on T}
 {}_{\sigma,\tau}T: \ \ \ r \cdot t \cdot r' : = \sigma(r)\tau(r') t = t^1r' \o r t^2.
\end{equation} 

Within $T$ is the (possibly zero) left ideal of Casimir elements
$\mathcal{C} := (A \o_B A)^A$ satisfying a type of right integral condition ($e,e' \in
\mathcal{C}$),
$$ee' = e \mu(e') = \mu(e') e$$
where $\mu: A \o_B A \to A$ denotes the multiplication map
$\mu(a \o a') := aa'$ an $A$-$B$ or $B$-$A$-split epimorphism. 
(Indeed, a right integral condition with respect to the $R$-bialgebroid 
structure on $T$ for a depth two extension $A \| B$.)
We note that the restriction $\mu: T \to R$ has image in the centralizer,
the restriction being denoted by the counit $\eps_T$ in section~2.

In section~4 we make use of two $R$-bimodules on $(A \o_B A)^A$
and $\hat{A}$. 
The Casimirs form an  $R$-$R$-bimodule ${}_R\mathcal{C}_R$
 as a submodule of
(\ref{inner structure on T}) given by
\begin{equation}
\label{eq: R-R-bimodule}
r \cdot e \cdot r' = e^1 r' \o r e^2.
\end{equation} The bimodule dual of $A$ forms
an $R$-$R$-bimodule ${}_R\hat{A}_R$ as a submodule
 of (\ref{eq: right R^e structure on S}).

We need a word about general modules as well.  Given any ring $\Upsilon$,
a right module $M_{\Upsilon}$ is isomorphic to a direct summand of another
right module $N_{\Upsilon}$ if we use the suggestive notation
$M_{\Upsilon} \oplus * \cong N_{\Upsilon}$.  A bimodule
${}_{\Upsilon}M_{\Upsilon}$ is of course the same as the left
or right $\Upsilon^e := \Upsilon \o_{\Z} \Upsilon^{\rm op}$-module $M$.   

Hirata extends Morita theory in an elegant way \cite{Hi} by defining
two modules $M_{\Upsilon}$ and $N_{\Upsilon}$ to be \textit{H-equivalent}
if both $M_{\Upsilon} \oplus * \cong N^n_{\Upsilon}$ (which equals
$N \oplus N \oplus \cdots \oplus N$, $n$ times) and
$N_{\Upsilon} \oplus * \cong M^m_{\Upsilon}$ for some positive integers
$n$ and $m$.  A  theory is outlined in \cite{Hi}, the main point
being that $\End M_{\Upsilon}$ and $\End N_{\Upsilon}$
are Morita equivalent rings with Morita context being
$\Hom (M_{\Upsilon},N_{\Upsilon})$ and $\Hom (N_{\Upsilon},M_{\Upsilon})$ with composition in either order
\cite{Lam}.


\section{Depth two theory}

Any ring extension $A \| B$ satisfies the property 
$A \oplus * \cong A \o_B A$ as either $B$-$A$ or $A$-$B$-bimodules,
since $\mu: A \o_B A \to A$ is an epi split by either $a \mapsto 1 \o a$
or $a \mapsto a \o 1$, respectively.  
Its converse from the point of view of Hirata's theory is satisfied
by special extensions introduced in \cite{KS} called depth two
extensions (due to their origins in subfactor theory).  Thus
a \textit{depth two extension} or D2 extension $A \| B$ is a ring extension
where 
\begin{equation}
\label{eq: D2}
A \o_B A \oplus * \cong A^n
\end{equation}
 as natural $B$-$A$-bimodules (\textit{left
D2}) and as $A$-$B$-bimodules (\textit{right D2}). Of course,
if some $A \| B$ satisfies the left and right conditions for different integers $m$ and
$n$, they both hold for the integer $\max \{ m,n \}$.  

There are several classes of examples of depth two extension $A \| B$;
such as an 
algebra that is finitely generated projective over its
base ring, finite Hopf-Galois extension,
H-separable extension or centrally projective extension \cite{KN,KS}.

How do we recognize depth two in any extension? This is not always easy from the definition
we provided above: let us provide a theorem with several equivalent conditions for left depth
two extensions (a similar theorem is then easily seen for right D2).
The theorem will use the left and right sides of the bimodules
${}_{\lambda,\, \rho}S$ and $T_{\sigma,\, \tau}$ respectively, given in
eqs.~(\ref{eq: left R^e structure on S}) and~(\ref{eq: outer structure on T}).

\begin{theorem}
\label{th-equiv}
The following are equivalent to the left depth two condition~(\ref{eq: D2})
 on a ring extension $A \| B$ in (1)-(4), on a Frobenius extension in (5)
and on a separable algebra $A$ with separable subalgebra $B$ over a field in (6):
\begin{enumerate}
\item the bimodules ${}_BA_A$ and ${}_BA \! \o_B \! A_A$ are $H$-equivalent;
\item there are $\{ \beta_j\}_{j=1}^n \subset S$ and $\{ t_j \}_{j=1}^n
 \subset T$ such that
\begin{equation}
\label{eq: quasibase}
a \o a' = \sum_j t_j \beta_j(a)a' \, ;
\end{equation}
\item as natural $B$-$A$-bimodules $A \o_B A \cong \Hom ({}_RS,{}_RA)$ and ${}_RS$ is a f.g.\
projective module;
\item as natural $B$-$A$-bimodules $T \o_R A \cong A \o_B A$  and $T_R$ is f.g.\ projective;
\item the endomorphism ring Frobenius extension has dual bases elements in $T$;
\item as a natural transformation between functors
from the category of right $B$-modules into the category
of right  $A$-modules, there is a natural monic from
$\Ind^A_B \Res^B_A \Ind^A_B$ into $N \Ind^A_B$
for some positive integer $N$. In particular, for
each pair of simple modules $V_B$ and $W_A$, the number
of isomorphic copies of $W$,  
\begin{equation}
\label{eq: ind} 
\bra \, \Ind^A_B \Res^B_A \Ind^A_B V , \, W \, \ket \leq
N\bra \, \Ind^A_B V, \, W \, \ket .
\end{equation}
\end{enumerate}
 \end{theorem}
\begin{proof}
We have already remarked above that the right $B^{\rm op} \o A$-modules $A$ and $A \o A$
are H-equivalent iff $A \| B$ is left D2. The second condition is called the D2 quasibase condition
and is noted in \cite{KS}.  It is based on condition~(\ref{eq: D2}) and the identifications
$\Hom ({}_B A_A, {}_BA \! \o \! A_A) \cong T$, $\Hom ({}_BA \! \o \! A_A, {}_BA_A)
\cong S$ and the existence of $2n$ maps $f_j$, $g_j$ in these two Hom-groups satisfying
$\sum_j f_j \circ g_j = \id_{A \o A}$.

Next we note the $B$-$A$ homomorphism $\chi: A \o A \stackrel{\cong}{\to} \Hom({}_RS, {}_RA)$ given by
$$\chi( a \o a')(\alpha) = \alpha(a)a'. $$
It is noted in \cite{KS} that   $\chi$ has inverse $\chi^{-1}(f) = \sum_j t_j f(\beta_j)$
if $A \| B$ is left D2.  For example, we check that $\chi \circ \chi^{-1} = \id$ on an $\alpha \in S$:

$$\sum_j \alpha(t_j^1) t_j^2 f(\beta_j) = f(\sum_j \alpha(t^1_j)t_j^2\beta_j) = f(\alpha)$$
for each $\alpha \in S$, since $\alpha(t_j^1)t_j^2 \in R$ for each $j$ and $a \o 1 = \sum_j t_j \beta_j(a)$ from Eq.~(\ref{eq: quasibase})
 implies $\alpha = \sum_j
\alpha(t_j^1)t_j^2 \beta_j$. We also note
that ${}_RS$ is finite projective with dual (projective) bases $\{ \psi(b_j) \}$, $\{ \beta_j \}$ 
where $\psi: T_R \stackrel{\cong}{\to} \Hom ({}_RS, {}_RR)$ is defined by $\psi(t)(\alpha) =
\alpha(t^1)t^2$. 

Conversely, if ${}_R S \oplus * \cong {}_RR^m$, then applying the functor $\Hom (-, {}_RA)$
from left $R$-modules into $B$-$A$-bimodules, we obtain from the isomorphism of the $A$-dual
of $S$ with the tensor-square that 
$A \o_B A \oplus * \cong {}_BA_A^m$.

The next condition follows from the D2 condition \cite{KS}, since the $B$-$A$-homomorphism
\begin{equation}
\label{eq: emm}
m: T \o_R A \stackrel{\cong}{\longrightarrow} A \o_B A, \ \ \  m(t \o a) = ta
\end{equation}
has an inverse $m^{-1}(a \o a') = \sum_j t_j \o \beta_j(a')a$ according to Eq.~(\ref{eq: quasibase}).
By the same token, $T_R$ is finite projective with dual bases
 $\{ t_j \}$, $\{ \bra \beta_j \,|\, -\ket \}$
from the nondegenerate pairing given in Eq.~(\ref{eq: pairing}) below. 

Conversely, if  $T_R \oplus * \cong R_R^t$ for some positive integer $t$, then applying
the functor $-\o_R A$ from right $R$-modules into $B$-$A$-bimodules results in 
${}_BA\! \o_B\! A_A \oplus * \cong {}_BA_A^t$, after applying the isomorphism of the tensor-square
with $T \o_R A$.  

The next-to-last condition depends on the fact that a Frobenius extension $A \| B$
with Frobenius homomorphism $E: A \to B$ and dual bases $x_i, y_i$
 has
Frobenius structure on its endomorphism ring extension  $\lambda:
A \into \End A_B \cong A \o_B A$
(isomorphism $f \mapsto \sum_i f(x_i) \o y_i$ with inverse
$a \o b \mapsto \lambda(a) E \lambda(b)$). Denote the endomorphism ring
structure induced on $A \o_B A$ by $A_1$; if $e_1 := 1 \o 1$,
its multiplication is the
$E$-multiplication given by
$$(a e_1 b)(c e_1 d) = aE(bc) e_1 d = a e_1 E(bc)d.$$  A Frobenius homomorphism
$E_A: A_1 \to A$ is given by the multiplication mapping $E_A(ae_1b) = ab$, 
since $x_i e_1 , e_1 y_i$ are dual bases for this mapping. Note that
$T = C_{A_1}(B)$, the centralizer of $B$ in the composite ring extension
$B \to A \into A_1$.    

The forward implication is \cite[Prop.\ 6.4]{KS}; viz., given a left D2
quasibase as above $t_i \in T, \beta_i \in S$ for $A \| B$, then
dual bases for $E_A$ in $T$ are given by $$\{ t_i \}, \ \\ \{ \sum_j \beta_i(x_j)  e_1 y_j \}$$
as shown in its proof.  (Also, $A \| B$ is necessarily right D2.)

For the reverse implication, we start with dual bases $\{ c_j \}$
and $\{ t_j \}$ in $T$ for $E_A$, and note that for every $a,b \in A$
\begin{eqnarray*}
ae_1 b & = &  \sum_j c^1_j e_1 c^2_j E_A(t_j a e_1 b) \\
       & = &  \sum_j c_j t^1_jE(t^2_ja)b \\
       & = & \sum_j c_j \eta_j(a)b 
\end{eqnarray*}
where $\eta_j = t^1_j E(t^2_j -)$ is clearly in $S$.  Thus, $\{ c_j \}$
and $\{ \eta_j \}$ form a left D2 quasibase.
(Similarly, a right D2 quasibase is given by $\{ t_j \}$
and $\{ \rho_j := E(- c^1_j)c^2_j \}$.)

The last condition for left depth two extension
of separable algebras follows in the forward implication
by tensoring the defining Eq.~(\ref{eq: D2}) by any
module $M_B$ from the left
and applying the definition of induced module:
$\Ind^A_B M := M \o_B A$,
and $\Res^B_A$ denoting usual restriction of an $A$-module to a  $B$-module.
The reverse implication follows from applying the functors and natural transformation to the regular
representation $B_B$.  Then 
$$B \o_B A \o_B A_A \into N B \o_B A_A = B \o_B A^N_A.$$
By naturality (applied to each morphism $\lambda_b: B_B \to
B_B$
for $b \in B$), this monic is also left $B$-linear.
Identify $B \o_B A \cong A$ as $B$-$A$-bimodules. Since
$B^{\rm op} \o A$ is a separable algebra, 
 Maschke's theorem shows that 
the $B$-$A$-monomorphism $A \o_B A \into
A^N$ splits.  Hence, $A$ is left D2 over $B$.
The Eq.~(\ref{eq: ind}) for the semisimple
category of right $A$-modules follows easily from this.  
\end{proof}

In addition to their ring structures, the constructs $S = \End {}_BA_B$ and $T= (A \o A)^B$ introduced in the preliminaries carry 
$R$-bialgebroid structures that are isomorphic to each others $R$-dual bialgebroids. (See \cite[2.4]{KS} for definition of left and right
$R$-dual bialgebroids over any ring $R$ given a finite projectivity condition). 
First, with respect to the  $R$-$R$-bimodule ${}_{\lambda,\, \rho} S$ in 
(\ref{eq: left R^e structure on S}), $S$ satisfies $$S \o_R S \cong \Hom ({}_B A\! \o\! A_B, 
{}_BA_B)$$ via 
$$\alpha \o \beta \longmapsto (a \o a' \mapsto \alpha(a) \beta(a')).$$
It follows that $S$ is an $R$-coring with coassociative comultiplication $\cop: S \to S \o_R S$
given by  
$$ \cop(\alpha)(a \o a') = \alpha(aa') $$
after identification, and counit $\eps_S: S \to R$ given
by $$ \eps_S(\alpha) = \alpha(1),$$ 
(see \cite{Brz} for the theory of corings). The ring-coring bimodule
structure $(S,R, \lambda, \rho, \cop, \eps)$ satisfies the axioms of a left bialgebroid in \cite{KS}
among which we find the important axioms $$\cop(\alpha \beta) = \cop(\alpha) \cop(\beta)$$
since $\Im \cop$ is contained in a subgroup of $S\o_R S$ where ordinary tensor product multiplication
makes sense; in addition, $\eps_S$ satisfies the following
alternative homomorphism property: $$\eps_S(\alpha \beta) = 
\eps_S(\alpha \lambda (\eps_S(\beta))).$$  

A \textit{left integral} $\ell$ in $S$ is an element in $S$ satisfying
\begin{equation}
\label{eq: left integral}
\alpha \circ \ell = \lambda(\eps_S(\alpha)) \ell
= \rho(\eps_S(\alpha)) \ell , \ \ \ (\forall \alpha \in S)
\end{equation}
For example, any element in $\hat{A}$ is a left integral.  A left integral
$\ell$ is \textit{normalized} if $\eps_S(\ell) = 1$, which for
$E \in \hat{A}$ is the case when $E(1) = 1$, i.e., $E: {}_BA_B \to {}_BB_B$
is a bimodule projection onto $B$, or \textit{conditional expectation} for
the \textit{split extension} $A \| B$.  Note the following
for any proper ring extension $A \| B$:

\begin{lemma}
 $A \| B$ is a split extension if and only if $\hat{A} \o_S A \cong B$
as $B$-$B$-bimodules and $\hat{A}_S$ is f.g.\ projective.  
\end{lemma}
\begin{proof}
($\Rightarrow $) Define $\hat{A} \o_S A \stackrel{\cong}{\to} B$
via $F \o a \mapsto F(a)$, with inverse $b \mapsto E \o b$.  
Note that the conditional expectation $E$ is a left identity.  

($\Leftarrow$)  If $\hat{A}_S \oplus * \cong S_S^n$, apply the
functor $- \otimes_S A_{B^e}$, obtaining
$$ {}_BB_B \oplus * \cong {}_BA_B^n . $$
Then there are $E_i \in \hat{A}$, $\lambda(r_i) \in \Hom ({}_BB_B, {}_BA_B) \cong
R$ such that $E := \sum_i E_i \circ \lambda(r_i)$ is $ \id_B$
when restricted to $B$.  
\end{proof}

Dually, we have $$S \cong \Hom (T_R, R_R)$$ via the nondegenerate pairing 
\begin{equation}
\label{eq: pairing}
\bra \alpha \, |\, t \ket =
\alpha(t^1)t^2
\end{equation} 
for depth two extension $A \| B$.  Consistent with this duality, $T$ has
the right bialgebroid structure $(T,R, \sigma, \tau, \cop_T, \eps_T)$ coming from
its ring structure and maps $\sigma, \tau:\, R \to T$ introduced in the preliminaries and
$R$-coring structure $(T, \cop_T, \eps_T)$ with respect to the bimodule $T_{\sigma,\, \tau}$
given in (\ref{eq: outer structure on T}), where comultiplication 
$$\cop(t) = t^1 \o 1_A \o t^2$$
under the identification $T \o_R T \cong (A \! \o \! A \! \o \! A)^B$ by $t \o t' \mapsto
t^1 \o t^2 {t'}^1 \o {t'}^2$, 
and counit $$\eps_T(t) = t^1 t^2.$$
The details may be found in \cite{KS}, but the bialgebroid structures are not very important 
in the present paper. 

A normalized right integral in $T$ is an element $u \in T$ such
that for every $t \in T$, $ut = u \sigma(\eps_T(t))= u \tau(\eps_T(t))$
with normalization condition $\eps_T(u) = 1_A$. 
We note the following easy proposition based on using the separability element
viewed in $T$ or a conditional expectation viewed in $\hat{A} \subset S$:
\begin{prop}
If $A \| B$ is a (resp., split) separable extension, then $T$ (resp., $S$) has a
normalized right (resp., left) integral.  
\end{prop}
 
Given a depth two extension $A \| B$, there are two measuring actions of $S$ and $T$. 
First, $T$ acts on $\mathcal{E} := \End {}_BA$
via the right action $f \ract t = t^1 f(t^2 - )$.  The invariants $\mathcal{E}^T$ under $\ract$
are defined to be $$\{ f \in \mathcal{E} \, | \, \forall \, t \in T, f \ract t = \lambda(\eps_T(t)) \circ f \}$$  
It is shown in \cite{KS} that $\mathcal{E}^T = \rho(A)$, the right multiplications by elements of 
$A$ (which is clear in one direction).  Forming the smash product $T \ltimes \mathcal{E}$,
one may show that it is isomorphic as rings 
to $\End {}_AA \! \o_B \! A$ \cite{K2002}. It is interesting to 
note that if $A$ additionally has a separability element $e \in T$,
then the mapping $\mathcal{E} \to A^{\rm op}$ defined by $f \mapsto f \ract e$
is a conditional expectation (a known fact recovered).  

Second, there is a very natural action of $S$ on $A$ given by $\alpha \lact a := \alpha(a)$.
Using a similar definition of invariants, we easily see that $A^S \supseteq B$.  
In addition to the depth two condition on $A \| B$, we
 need an extra condition on $A_B$ to prove that $B = A^S$:  in \cite{KS}
this is done by assuming $A_B$ a balanced module, i.e., $B \cong \End {}_{\mathcal{E}'}A$
where $\mathcal{E}' = \End A_B$ and we consider only natural actions.  
(For example, if $A_B$ is a generator, it is balanced.)
The smash product $A \rtimes S$ is isomorphic to $\mathcal{E}'$ \cite{KS}.
 
We end this section with a necessary condition for an extension $A \| B$ to be left D2.
Recall the left $T$-module $A \o_B A$ and the module $R_T$ both defined in the preliminaries.
We let $Z(A)$ denote the center subring of any ring $A$.  

\begin{theorem}
\label{th-square}
Suppose $A \| B$ is D2.  Then $R \o_T (A \o_B A) \cong A$ as $A$-$A$-bimodules
and $\End R_T \cong Z(A)$.
\end{theorem}
\begin{proof}
The mapping $\gamma: R \o_T (A \o_B A) \to A$ given by
\begin{equation}
\label{eq: gamma}
\gamma(r \o a \o a') = ara'
\end{equation}
is an $A$-$A$-bimodule epimorphism for every ring extension $A \| B$
(a type of ternary product through which the multiplication
mapping $\mu: A \o_B A \rightarrow A$ factors).

We work with the bimodule ${}_TT_R$ given by $$t \cdot t' \cdot r := tt' \tau(r) = {t'}^1 t^1 \o t^2 {t'}^2 r$$
derived from the structures introduced in the preliminaries. Note that
the diagram on this page below  is a commutative square, where the left vertical and bottom
horizontal arrows are both induced from canonical maps of the type $R \o_R A \cong A$,
$r \o a \mapsto ra$, while the top horizontal mapping is induced from $m: T \o_R A \cong A \o A$
given in Eq.~(\ref{eq: emm}). Both sides of the square send $r \o t \o a$ in the upper left-hand
corner into $t^1rt^2a \in A$.  It follows that the fourth map $\gamma$ is an isomorphism.

$$\begin{diagram}
R \o_T T \o_R A &&\rTo^{\cong}_{\id \o m} && R \o_T (A \o_B A)\\
\dTo^{\cong} && && \dTo_{\gamma}\\
R \o_R A && \rTo^{\cong} && A 
\end{diagram}$$   

\vspace{.5cm}

It follows that $\gamma^{-1}(a) = 1 \o_T ( 1 \o a)$, a fact that is directly checked as follows:
$$ 1 \o_T (1 \o ara') = 1 \o (1 \o \sum_j t_j^1 r t_j^2 \beta_j(a)a') = \sum_j t_j^1 r t_j^2 \o_T 
( 1 \o \beta_j(a)a') $$
$$ = \sum_j r \o_T ( t_j \beta_j(a)a') = r \o (a \o a').$$

Clearly, $Z(A) \into \End R_T$ via $\lambda$.  Conversely, if $f \in \End R_T$,
we note that $f(1)$ satisfies a type of left integral condition, $$ t^1 t^2 f(1) =
f(1) \cdot (t^1 t^2 \o 1) = f(1 \cdot t) = f(1) \cdot t = t^1 f(1) t^2. $$

Next we recall the right action of $T$ on $\mathcal{E} = \End {}_BA$ outlined above:
we see that $R_T \into \mathcal{E}_T$ is a submodule via $\lambda$ since
for every $t \in T$, 
$$\lambda(r) \ract t = \lambda(t^1 r t^2) = \lambda(r \cdot t) . $$
Then $f(1) \in R$ is in the invariants $\mathcal{E}^T = \rho(A)$, since
$$ \lambda(f(1)) \ract t = \lambda(t^1 f(1) t^2) = \lambda(\eps_T(t)) \lambda(f(1)),$$
whence $\lambda(f(1)) = \rho(a) \in \mathcal{E}^T$. Now $a = f(1)$
by evaluating at $1$, and $a \in Z(A)$ by evaluating at any $b \in A$. 
\end{proof}

\begin{example}
\begin{sl}
We have identified examples of depth two extensions coming from centrally projective, H-separable or finite
Hopf-Galois extensions of rings or algebras.  
Let us see by means of the theorem  an easy example of an algebra 
extension which is not D2.  If $A$ is the algebra of upper triangular $2 \times 2$ matrices
over any field, $B$ the subalgebra of diagonal matrices, then $R = B$,
$T$ is spanned by $\{ e_{11} \otimes e_{11},\, e_{22} \otimes e_{22} \}$
in terms of matrix units $e_{ij}$ and 
it is easy to see that $$\dim \End (R_T) = \dim \End R_R = \dim R = 2$$
while $\dim Z(A) = 1$. 
There are also nonexamples of depth two extensions coming
from various natural subalgebras of finite dimensional exterior algebras (where condition (4) in Theorem~\ref{th-equiv} is useful as well as the fact that
split depth two extensions are f.g.\ projective by
eq.~\ref{eq: quasibase}, therefore free over local rings).   
\end{sl}
\end{example} 


\section{Depth two as normality for subgroups and Hopf subalgebras}

In this section we show that depth two subgroup algebras and Hopf subalgebras are closely
related to normal subgroups and normal Hopf subalgebras.  For subgroups of finite groups over
the complex numbers we have a definitive result in Theorem~\ref{th-burkhard} and its corollary. 

\subsection*{Character theory \cite{Alp}}
Let $H$ be a subgroup of finite index in a group $G$, and let $R$ be a
commutative ring. We are interested in the question when the extension
of group algebras $RG |RH$ is D2. First, we recall that $RG | RH$ is D2
whenever $H$ is normal in $G$ \cite[3.9]{KS}. Indeed, suppose that $\{g_1,
\ldots, g_n\}$ is a transversal for the cosets of $H$ in $G$. Then the
map
$$(RG)^n \longrightarrow RG \otimes_{RH} RG, \quad (x_1, \ldots, x_n)
\longmapsto \sum_{j=1}^n x_jg_j^{-1} \otimes g_j,$$
is an isomorphism of $RG$-$RH$-bimodules. Similarly, $(RG)^n$ and $RG
\otimes_{RH} RG$ are isomorphic $RH$-$RG$-bimodules.


For the converse, let $R$ be the field of complex numbers $\C$.
Suppose that $\chi_1,\ldots,\chi_n$ are the irreducible characters of 
a finite group $G$,
while $\psi_1,\ldots,\psi_m$ are the irreducible characters of a subgroup $H$.
Suppose $a^i_j$ are nonnegative integers in an induction-restriction table
for $H \leq G$ such that $$\psi_i^G = \sum_{j=1}^n a^i_j \chi_j,$$
and $c^r_s$ are nonnegative integers satisfying
$$(((\psi_r)^G)_H)^G = \sum_{s=1}^n c^r_s \chi_s. $$
Then the group algebra extension $A \| B$  is D2 if and only if 
there is a positive integer $N$ such that $c^i_j \leq Na^i_j$ for
all $i= 1,\ldots,m$ and $j= 1, \ldots,n$.

In particular, $A \| B$ is not D2 if $a^i_j = 0 $ while $c^i_j \neq 0$.
For example, the induction-restriction table (based
on Frobenius reciprocity $(\psi_i^G, \chi_j)_G = (\psi_i,(\chi_j)_H)_H$) for the permutation
groups $S_2 \leq S_3$ is given by
$$
\begin{tabular}{l|lll} 
$S_2 \leq S_3$ & $\chi_1$ & $\chi_2$ & $\chi_3$ \\ \hline
$\psi_1$  & 1 & 0 & 1 \\
$\psi_2$ & 0 & 1 & 1 
\end{tabular}
$$
where $\psi_1$, $\chi_1$ denote the trivial characters, $\psi_2$, $\chi_2$
the sign character, and $\chi_3$ the two-dimensional 
irreducible character of $S_3$.
It follows that $\psi_1^G = \chi_1 + \chi_3$, while 
$$ ((\psi_1^G)_H)^G = 2\psi_1^G + \psi_2^G = 2 \chi_1 + \chi_2 + 3 \chi_3.$$
Hence,  $ ((\phi_1^G)_H)^G$ is not even a subcharacter of $\psi_1^G$ or any
of its integral multiples.  

The last example and results from operator algebras raise the question
whether subgroups of finite groups that are depth two as complex group subalgebras
are necessarily normal subgroups.  The next theorem answers this affirmatively.  

\begin{theorem}
\label{th-burkhard}
 Let $H$ be a subgroup of a finite group $G$ such that      
$\C G \| \C H$ is a depth two ring extension. Then $H$ is normal in     
$G$.                                                                      
\end{theorem}                                                              
\begin{proof}
 The depth two hypothesis implies that there exists a positive     
integer $n$ such that                                                          
$$\langle {\rm Ind}_H^G ({\rm Res}_H^G ({\rm Ind}_H^G (\psi))) \|               
\chi \rangle_G \le n \langle {\rm Ind}_H^G (\psi) \| \chi \rangle_G$$           
for $\chi \in {\rm Irr}(G)$ and $\psi \in {\rm Irr}(H)$. 
We choose $\chi = 1_G$ and $1_H \ne \psi \in {\rm Irr}(H)$. Then, by    
Frobenius reciprocity, we have                                                                                              
$$\langle {\rm Ind}_H^G (\psi) \| 1_G \rangle_G = \langle \psi \| 1_H \rangle_H  
= 0.$$                                                                         
Hence, we have                                                                 
$$0 = \langle {\rm Ind}_H^G ({\rm Res}_H^G ({\rm Ind}_H^G (\psi))) \| 
1_G \rangle_G = \langle {\rm Res}_H^G ({\rm Ind}_H^G (\psi)) \| 1_H \rangle_H  $$
$$= \sum_{HgH \in H \backslash G/H} \langle {\rm Ind}_{H \cap gHg^{-1}}^H       
({\rm Res}_{H \cap gHg^{-1}}^{gHg^{-1}} ({^g\psi})) \| 1_H \rangle_H, $$     
by Frobenius reciprocity and Mackey's formula. We conclude that, for $g \in G$,
\begin{eqnarray*}
0 & = & \langle {\rm Ind}_{H \cap gHg^{-1}}^H ({\rm Res}_{H \cap       
gHg^{-1}}^{gHg^{-1}} ({^g\psi})) \| 1_H \rangle_H  \\                             
&=& \langle {\rm Res}_{H \cap gHg^{-1}}^{gHg^{-1}} ({^g\psi}) \|1_{H \cap         
gHg^{-1}} \rangle_{H \cap gHg^{-1}} \\
& =& \langle {\rm Res}_{g^{-1}Hg \cap H}^H (\psi) \| 1_{g^{-1}Hg \cap H}          
\rangle_{g^{-1}Hg \cap H} \\                                                     
& = & \langle \psi \| {\rm Ind}_{g^{-1}Hg \cap H}^H (1_{g^{-1}Hg \cap H})           
\rangle_H, 
\end{eqnarray*}                                                              
again by Frobenius reciprocity and conjugation. On the other hand, we have     
$$\langle 1_H \| {\rm Ind}_{g^{-1}Hg \cap H}^H (1_{g^{-1}Hg \cap H})
\rangle_H = \langle 1_{g^{-1}Hg \cap H} \| 1_{g^{-1}Hg \cap H} \rangle_{g^{-1} 
Hg \cap H} = 1,$$                                                             
using Frobenius reciprocity one more time. Thus                               
$${\rm Ind}_{g^{-1}Hg \cap H}^H (1_{g^{-1}Hg \cap H} ) = 1_H.$$               
Comparing degrees we get $H = g^{-1}Hg \cap H$. Hence $H$ is normal in $G$.      
\end{proof} 

We summarize the theorem with its converse proven above using characters (and group-theoretically
in \cite{KS}
for finite index normal subgroups over general fields).

\begin{cor}
Let $H$ be a subgroup of a finite group $G$.  Then $ \C G \| \C H$ is a depth two ring extension if
and only if $H$ is a normal subgroup of $G$.
\end{cor}  

It is well-known that a finite group algebra is a Hopf-Galois extension of a normal subgroup algebra:
see our next subsection for an exposition of this fact in a more general setting.  
Theorem~\ref{th-burkhard} implies a converse, since a finite Hopf-Galois extension is a depth two extension.  

\begin{cor}
Let $H$ be a subgroup of a finite group $G$.  If $ \C G \| \C H$ is a Hopf-Galois extension, then 
 $H$ is a normal subgroup of $G$.
\end{cor}                  

\subsection*{Normal Hopf subalgebras \cite{Mo} are depth two} 
Consider any finite-dimensional Hopf algebra $H$ over a field $k$
 with antipode $S: H \to H$, counit $\eps: H \to k$
and comultiplication notation $\cop(a) = a\1 \o a\2$ for each $a \in H$.  A Hopf subalgebra $K$ is \textit{normal} if $S(a\1)xa\2
\in K$ and $a\1 x S(a\2) \in K$ for each $a \in H$ and $x \in K$.  Form the subset $K^+ := \ker \eps \cap K$
and the left ideal $HK^+$ in $H$.  Note that $HK^+ = K^+ H$ since given $x \in K^+$ and $a \in H$
$$xa = \eps(a\1)xa\2 = a\1 S(a\2) x a\3 \in HK^+,$$
the converse $HK^+ \subseteq K^+H$ being shown similarly.
It follows that $HK^+$ is a Hopf ideal and $\overline{H} := H/ HK^+$ is a Hopf algebra.  There is the natural surjective Hopf algebra homomorphism $H \to \overline{H}$ denoted
by $a \mapsto \overline{a}$. 

\begin{prop}
Let $H$ be a finite-dimensional Hopf algebra and $K$ a normal Hopf subalgebra in $H$.  Then $H \| K$ is an
$\overline{H}$-Galois extension, therefore of depth two.
\end{prop}
\begin{proof}  The natural coaction $H \to H \o_k \overline{H}$
given by $a \mapsto a\1 \o \overline{a\2}$ makes $H$ into an $\overline{H}$-comodule algebra with coinvariants  $K$  \cite[3.4.3]{Mo}.  The $\overline{H}$-extension $A \| B$ is
Galois \cite[p.\ 30]{Mo}
(with Galois mapping $\beta: H \o_K H \to H \o \overline{H}$ given by $\beta(x \o y) = xy\1 \o \overline{y\2}$).

Finally recall from \cite{KS, K2003} that any finite Hopf-Galois extension $A \| B$ is a
 depth two extension.
\end{proof}

For example, let $G$ be a finite group with normal subgroup $N$.  Then $H = kG$ and $K = kN$ is a normal Hopf
subalgebra pair with Galois mapping $\beta : H \o_K H \to H \o k[G/N]$ given by $\beta(g \o g') = gg' \o g'N$
for $g, g' \in G$.  

We pose a problem of whether a converse of the proposition holds; e.g. if a depth two Hopf subalgebra or
 depth two semisimple Hopf subalgebra pair is normal.


\section{A fresh look at  separability}

In this section we meet a new criterion for separability which is very close to
the necessary condition in Theorem~\ref{th-square} for depth two extensions.  We will 
work with a \textit{separable extension} $A \| B$ which is defined
to be a ring extension where $\mu: A \o_B A \to A$ is split as an
$A$-$A$-epimorphism.  The image of $1_A$ under any splitting is
called a separability element $e \in A \o_B A$ which of course
satisfies $\mu(e) = 1$ and $e \in (A \o_B A)^A$. Fix this (nonunique) separability
element $e = e^1 \o e^2$ (suppressing a possible summation) for any separable extension below.   

\begin{theorem}
\label{th-newsep}
A ring extension $A \| B$ is separable if and only if
$R_T$ is f.g.\ projective and  $R \o_T (A \o_B A) \cong A$
as natural $A$-$A$-bimodules.
\end{theorem}
\begin{proof}
($\Rightarrow$) An inverse to the $A$-$A$-homomorphism 
$\gamma$ defined in Eq.~(\ref{eq: gamma})
is given by $\gamma^{-1}(a) = 1 \o ea \in R \o_T (A \o_B A)$. We compute
 $\gamma \circ \gamma^{-1} = \id$
from $e^1 e^2 a = a$, and $\gamma^{-1} \circ \gamma = \id$ from 
$$1 \o_T eara' = 1 \o_T er \cdot (a \o a') = 1 \cdot er \o (a \o a') = r \o (a \o a'). $$

$R_T$ is finite projective $\Leftrightarrow$ there are finite dual bases
$\{ r_i \in R \}$, $\{ f_i \in \Hom (R_T, T_T) \cong T'\}$ via $f \mapsto f(1)$ and
where the subring $T'$ are $R$-centralizing left integral-type elements, i.e.,  
\begin{equation}
\label{eq: areprime}
T' = \{ p \in T | \, p\cdot t = t^1pt^2 = pt^1t^2 = t^1t^2p,\, \forall t \in T \} 
\end{equation}
$\Leftrightarrow$ there are $\{ p_i \in T' \}$, $ \{ r_i \in R \}$ such
that 
\begin{equation}
\label{eq: R_T projective}
r = \sum_{i=1}^n p_i^1 r_i p_i^2 r = \sum_i r p_i^1 r_i p_i^2 \ \ \ (r \in R)
\end{equation}
But $n =1$, $r_1 = 1$ and $p_1 = e$ will do.

($\Leftarrow$)  If $R_T \oplus * \cong T_T^m$, then tensoring by $- \o_T (A \o_B A)$,
we obtain ${}_AA_A \oplus * \cong {}_AA \! \o_B \! A_A^m$.  
It follows by a simple argument using split epimorphism, canonical maps between a
product and its components, and the 
identifications $\Hom ({}_AA_A, {}_A A \! \o_B \! A_A) \cong
(A \o_B A)^A$ and $\Hom ({}_A A \! \o_B \! A_A, {}_AA_A) \cong R$
(via $f \mapsto f(1 \o 1)$) that there are $m$ Casimir elements $e_i$
and $m$ elements $r_i \in R$ such that $\sum_{i=1}^m r_i \cdot e_i = 1_A$.
Then $e := \sum_{i=1}^m e_i^1 \o r_i e_i^2$ is a separability element.
\end{proof}

We next recall a special case of separable extension which has been extensively studied
by Ikehata,  Sugano, Szeto and others.
A ring extension $A \| B$ is said to be \textit{H-separable} (after Hirata \cite{Hi})
if  $A \o_B A \oplus * \cong A^n$ as $A^e$-modules.  Equivalently, there are
$n$ Casimir elements $e_i$ and $n$ elements $r_i \in R$, called
an \textit{H-separability system}, such that $1 \o 1 = \sum_{i=1}^n
r_i e_i$.  Hirata \cite{Hi} shows that $R$ is f.g.\ projective over $ Z := Z(A)$, therefore 
$R_{Z}$ is a generator module, 
\begin{equation}
\label{eq: isomorphism}
A \o_{Z} A^{\rm op}  \cong \Hom (R_{Z}, A_{Z})
\end{equation}
 via $a \o a' \mapsto \lambda(a)\rho(a')$, 
and so a separability element for $A \| B$ corresponds to a splitting map for $Z \into
R$, which exists by the generator property. It follows that $A$ and $A \o_B A$ are
H-equivalent $A$-$A$-bimodules, another characterization of an H-separable extension $A \| B$. 

\begin{theorem}
A ring extension $A \| B$ is H-separable if and only if $R_T$ is a generator and
$R \o_T (A \o_BA) \cong A$ as natural $A$-$A$-bimodules.
\end{theorem}
\begin{proof}
($\Rightarrow$) Since $A \| B$ is also a separable extension, $R \o_T (A \o_BA) \cong A$
by Theorem~\ref{th-newsep}.  Also $R_T$ is f.g.\ projective, but we need to see that
it is a generator (of the category of $T$-modules).
But $R_T$ generator $\Leftrightarrow$ there are finite $\{ r_i \in R \}$, $\{ g_i \in \Hom (R_T,T_T)
\cong T'\}$ where $T'$ is given by (\ref{eq: areprime})
such that 
\begin{equation}
\sum_i g_i(r_i) = 1_T \ \ \Leftrightarrow \ \  \sum_i q_i r_i = r_i q_i = 1 \o 1 \ \ (q_i = g_i(1)) 
\end{equation}
which  is satisfied by an H-separability system $q_i = e_i$, $\{ r_i \}$ as defined above.

($\Leftarrow$) If $T_T \oplus * \cong R_T^n$, which is the generator condition in one of its
equivalent forms \cite{Lam}, then tensoring by $- \o_T (A \o_B A)$ and applying the isomorphism,
we obtain ${}_A A\! \o_B \! A_A \oplus * \cong {}_AA_A^n$, which shows $A \| B$ is H-separable.
\end{proof}
\begin{cor}
If $A \| B$ is H-separable, then ${}_{Z}R_T$ is a faithfully balanced bimodule
with $Z$ and $T$ Morita equivalent rings.
\end{cor}
\begin{proof}
 By restricting $A$-modules to $B$-modules (or pullback in case of ring
arrow $B \to A$) we see from the defining conditions for H-separability and depth two
 that an H-separable extension $A \| B$ is depth two. Then $\End R_T \cong Z$ from
Theorem~\ref{th-square} and $\End {}_{Z}R \cong T$ from taking the 
$B$-centralizer of (\ref{eq: isomorphism}).  
Note that we have already seen in this section that ${}_{Z}R$ and $R_T$ are progenerators.
It follows that the rings $T$ and $Z$ are Morita equivalent.

 We may alternatively
see this last point
 from the fact noted above that $A$ and $A \o_B A$ are H-equivalent $A^e$-modules,
whence their endomorphism rings are Morita equivalent. Hirata theory also
shows handily that the Morita context bimodules here are $R$ as before, 
$(A \o_B A)^A$ with $T$-module as left ideal in $T$ and trivial $Z$-module, 
and associative pairings
 $\bra e \| r \ket = er \in T$, $[ r \| e ] = e^1 r e^2 \in Z$.
\end{proof}

\begin{example}
\begin{sl}
Suppose $A \| B$ is an H-separable right $K$-Galois extension for some Hopf algebra $K$.
Then $T^{\rm op} \cong R \rtimes K^{\rm op}$ as algebras, where
$K$ acts on the centralizer $R$ via the Miyashita-Ulbrich action \cite{K2003}.  Then the
 module $R_T \cong {}_{R \rtimes K^{\rm op}} R$ is a generator.  But this module
is a generator iff $R$ is a right $K^{\rm op\, *}$-Galois extension over
$R^K = Z$ by \cite[p.\ 133]{Mo}, which shows \cite[3.1]{K2001} by other means. 
\end{sl}
\end{example}

\begin{example}
\begin{sl}
A result of Sugano states that given an algebra $\Lambda$ over a commutative base
ring $R$ with center $C = Z(\Lambda)$, then $\Lambda$ is H-separable over $R$
iff $\Lambda$ is separable over $C$ and $C \o_R C \cong C$ via $c_1 \o c_2 \mapsto c_1 c_2$.  In other words, Azumaya algebras are H-separable over their
centers, while H-separable algebras are Azumaya algebras with center not much
larger than the base ring.  A depth two, separable extension $A \| B$ that
is not H-separable is then given by letting $K$ be any field, the direct product
algebra $A = K^n$
for integer $n > 1$, $B = K 1_A$.   
\end{sl}
\end{example}


\section{A trace ideal condition for Frobenius extensions}

Recall that a ring extension $A \| B$ is \textit{Frobenius} if ${}_BA_A \cong {}_B \Hom (A_B, B_B)_A$
(via a \textit{Frobenius isomorphism}, say $\psi$) 
and $A_B$ is f.g.\ projective.  Thus the right $B$-dual $A^*$ of $A$ is a free rank one right $A$-module.
A free generator $E: A \to B$ is called a \textit{Frobenius homomorphism}, and $E \in \hat{A}$ since
it is an image of an invertible element in $R$ under $\psi$: of course, we take 
$E = \psi(1)$.
Projective bases $\{ x_i\in A \}_{i=1}^n$, $\{ f_i\in A^* \}_{i=1}^n$ for $A_B$ convert to
dual bases $\{ x_i \}$, $\{ y_i := \psi^{-1}(f_i) \}$ satisfying
\begin{equation}
\label{eq: db}
\sum_{i=1}^n x_i E(y_i a) = a = \sum_{i=1}^n E(a x_i)y_i.
\end{equation}  
The first equation follows immediately from the dual bases equation for $A_B$, while
the second will follow from the first and injectivity of $\psi$ by
 evaluating $\psi(a - \sum_i E(ax_i)y_i)$ on any $x \in A$.
From the equations above, one computes that $\sum_{i=1}^n x_i \o y_i$ is a Casimir element:
$$ \sum_i ax_i \o y_i = \sum_{i,j} x_j \o E(y_j ax_i)y_i = \sum_j x_j \o y_j a.$$

There are various questions dating back to Nakayama and Eilenberg
of the type, to what extent separable extensions are Frobenius \cite{CK}?
The most important theorem in this area to date is the one of Endo-Watanabe
establishing that a separable algebra $A$ over a commutative ring $k$ is
a symmetric algebra (in particular, Frobenius algebra) if $A$ is faithful,
f.g.\ projective as $k$-module.  This theorem has been brought to bear on the
centralizer by Sugano to show that
various H-separable and centrally projective separable extensions are Frobenius extensions.

The next proposition gives necessary conditions for an extension to be Frobenius,
properties that could rule out various separable extensions from being Frobenius.
\begin{prop}
\label{prop-Frob}
Suppose $A \| B$ is a Frobenius extension.  Then $ \hat{A}_R$
and $(A \o_B A)^A$ are both free right $R$-modules of rank one.
\end{prop}
\begin{proof}
First  apply the functor $(-)^B$, the $B$-centralizer to both sides of $B$-$A$-isomorphism
$A \cong A^*$, obtaining $R\cong \hat{A}$ (since the  $B$-$A$-bimodule structure on $A^*$
is given by $b \cdot f \cdot a = \lambda(b) \circ f \circ \lambda(a)$,
so $f \in (A^*)^B \Leftrightarrow \lambda(b) f= f \lambda(b)$ for each $b \in B$).
This isomorphism of right $R$-modules (cf.\ preliminaries and module
 structure~(\ref{eq: right R^e structure on S}))
 is given by $r \mapsto E \cdot r$,
with inverse $F \mapsto \sum_i F(x_i)y_i$.  

We next note that 
\begin{equation}
A \o_B A \stackrel{\cong}{\longrightarrow} \End A_B, \ \ a \o a' \mapsto \lambda(a) \circ E \circ \lambda(a')
\end{equation}
Its inverse is given by $f \mapsto \sum_i f(x_i) \o y_i $ where $E, x_i, y_i$ are defined above.  

Taking the $A$-centralizer of this $A$-$A$-isomorphism, we see that 
\begin{equation}
(A \o_B A)^A \cong \End {}_AA_B \cong R^{\rm op}
\end{equation}
where the composed isomorphism comes out as
 $e \mapsto e^1E(e^2)$ with inverse $r \mapsto \sum_i x_ir \o y_i$. 
These maps are right $R$-module isomorphisms with respect to the natural right $R$-module 
and module~(\ref{eq: R-R-bimodule}). 
\end{proof}

\begin{example}
\begin{sl}
The full $n \times n$ matrix algebra $A$ over a field is a separable algebra, therefore
Frobenius algebra, and so $(A \o A)^A \cong A$.  Indeed, $\sum_i e_{ij} \o e_{ri}$
where $j,r = 1, \ldots,n$ form a basis of $n^2$ Casimir elements. ($n$ of these
are separability elements in $\mu^{-1}(1)$, and their average, if the characteristic
does not divide $n$, is the unique symmetric separability element.)
\end{sl}
\end{example}

We recall from section~1 that
 $\mathcal{C} := (A \o_B A)^A$ denotes the Casimir elements of $A \| B$, viewed below as a submodule of the $R$-$R$-bimodule
${}_{\sigma,\, \tau}T$; also recall the $B$-valued bimodule
homomorphisms $\hat{A}$, viewed below as a submodule  of
$S_{\lambda,\, \rho}$. 
Next define two $R$-$R$-homomorphisms on the tensor product of the $R$-$R$-bimodules
$\hat{A}$ and $\mathcal{C}$ in either order:
\begin{equation}
\label{eq: Psi}
\Psi:\ \hat{A} \o_R \mathcal{C} \to R, \ \ \ \Psi(F \o e) = e^1 F(e^2),
\end{equation}
 a well-defined $R$-$R$-homomorphism since $r \cdot F \cdot r' = F \circ \lambda(r')\circ \rho(r)$
($F \in \hat{A}$) and $r \cdot e \cdot r' = e^1 r' \o r e^2$ 
($e \in (A \o A)^A, r,r' \in R$), whence
$$\Psi(rF \o er') = e^1 r' F(e^2 r) = r e^1 F(e^2)r'.$$

Similarly,
\begin{equation}
\label{eq: Phi}
\Phi:\ \mathcal{C} \o_R \hat{A} \longrightarrow R, \ \ \ \Phi(e \o F) = F(e^1)e^2
\end{equation}
defines an $R$-$R$-homomorphism.

We first see that a necessary condition for $\Phi$ and $\Psi$ to be surjective
is that $\hat{A}$ and $\mathcal{C}$ be left and right $R$-generators (by the trace
ideal characterization of generator).

\begin{theorem}
The datum $(R,\, R,\, \hat{A},\, \mathcal{C},\, \Phi,\, \Psi)$ is a Morita context \cite{Lam}
where $\Phi$ and $\Psi$ are surjective if any one of the conditions below are satisfied
\begin{enumerate}
\item $A \| B$ is a Frobenius extension;
\item $R$ is a simple ring and either $\Phi \neq 0$ or  $\Psi \neq 0$;
\item one of $\hat{A}$ and $\mathcal{C}$ is an $R$-progenerator while the other is
isomorphic to its corresponding $R$-dual with $\Psi$ and $\Phi$ corresponding to evaluation and
co-evaluation.
\end{enumerate}
\end{theorem}
\begin{proof}
We first show that the two associativity squares below are commutative. First, the commutativity of the square, 

$$\begin{diagram}
\mathcal{C} \o_R \hat{A} \o_R \mathcal{C}&&\rTo^{1 \o \Psi} && \mathcal{C} \o_R R \\
\dTo^{\Phi \o 1} && && \dTo_{\rm can}\\
R \o_R \mathcal{C} && \rTo_{\rm can} && \mathcal{C}  
\end{diagram}$$ 
corresponds to the equality ($e,f \in (A \o A)^A;\, E,F \in \hat{A}$)
$$ e \cdot \Psi(F \o f) = e \cdot (f^1 F(f^2)) = e^1f^1F(f^2) \o e^2 $$
$$ = f^1 \o F(f^2 e^1)e^2 = (F(e^1)e^2)\cdot f =  \Phi(e \o F)\cdot f.$$
The commutativity of the other square, 
$$\begin{diagram}
 \hat{A} \o_R \mathcal{C}\o_R \hat{A} &&\rTo^{1 \o \Phi} && \hat{A} \o_R R \\
\dTo^{\Psi \o 1} && && \dTo_{\rm can}\\
R \o_R \hat{A} && \rTo_{\rm can} && \hat{A} 
\end{diagram}$$ 
follows from: 
$$ E \cdot \Phi(e \o F) = F(e^1) E(e^2 -) = F(-e^1) E(e^2) = \Psi(E \o e) \cdot F.$$

If $A \| B$ is a Frobenius extension, then there is $E \in \hat{A}$
and $e \in \mathcal{C}$ such that $\Phi(e \o E) = 1 = \Psi(E \o e)$,
as outlined above.  
But $\Im \Psi$ and $\Im \Phi$ are both two-sided ideals in $R$, since
both maps are $R$-$R$-homomorphisms. 

By the same token, if either map is nonzero and $R$ is simple,
the map is surjective.  The following is an elementary fact from
\cite{KP}: $$\Phi(e \o \hat{A}) \neq 0 \Leftrightarrow 
\Psi(\hat{A} \o e) \neq 0 $$
and similarly
$$ \Psi(E \o \mathcal{C}) \neq 0 \Leftrightarrow \Phi(\mathcal{C} \o E)
\neq 0. $$

The last item is a  consequence of  the well-known theorem of Morita
that if $P_R$ is a progenerator, then $R$ and $\End P_R$ are Morita
equivalent with context $P$ and its $R$-dual $P^*$ with evaluation
and coevaluation as the associative bimodule isomorphisms.  
\end{proof}

If we consider Picard groups for noncommutative rings as the group of auto-equivalences
of the corresponding full module category, we see via Morita theory
 that  $\Phi$ and
$\Psi$ are surjective iff  $\hat{A}$ and $(A \o A)^A$ represent inverse elements in the
Picard group of $R$, which coalesce by Proposition~\ref{prop-Frob} 
to the neutral element in case $A \| B$ is Frobenius. 

The theorem provides a criterion for ring extensions to be Frobenius in analogy
with the trace ideal condition for  modules to be  generators.  The images of
$\Psi$ or $\Phi$ are two-sided ideals, which we check to see if the maps are
surjective; if so, we can apply the next Frobenius result to right f.g.\ projective
extensions with centralizer $R$
satisfying $xy = 1$ for $x,y \in R$ implies $yx = 1$ (i.e., $R$ is \textit{Dedekind-finite}
\cite{Lam}). 
\begin{theorem}
Suppose $A \| B$ is a ring extension where $A_B$ is f.g.\ projective and its centralizer
is Dedekind-finite. 
If the $R$-bimodule pairings $\Psi$ and $\Phi$ on $\hat{A}$ and $\mathcal{C}$
 defined above are surjective, then $A$ is a Frobenius extension of $B$.
\end{theorem}
\begin{proof}
If the pairings above are surjective, there are $e,f \in \mathcal{C}$ and $E,F \in \hat{A}$
such that $$E(e^1) e^2 = 1 = f^1 F(f^2).$$ Then the homomorphism ${}_BA_A \to {}_BA^*_A$
given by $a \mapsto Ea$ (where $Ea(x) = E(ax)$ as usual) is injective, for
$$ Ea = 0 \Rightarrow 0 = (Ea)(e^1) e^2 = E(e^1)e^2a = a.$$
But the homomorphism ${}_BA_A \to {}_B A^*_A$ given by $a \mapsto Fa$ is onto
since
$$ a = af^1 F(f^2) = f^1 F(f^2 a) \ \Rightarrow \ \forall \phi \in A^*: \ \phi(a)= \phi(f^1)F(f^2 a) =
Fx(a) $$
where $x = \phi(f^1)f^2$.  It follows that  $b = E(f^1)f^2 \in R$ satisfies $Fb = E$. Then
$$ 1 = E(e^1)e^2 = F(be^1)e^2 = F(e^1)e^2 b $$
implies $b$ is invertible by Dedekind-finiteness of $R$.  
Hence our map $a \mapsto Ea$ is 1-1 and onto, since $x \mapsto Fx$ is onto and $Fx = Eb^{-1} x$ for
any $x \in A$. It follows that ${}_BA_A \cong {}_B \Hom (A_B, B_B)_A$ and $A_B$ is finite
projective, whence $A \| B$ is a Frobenius extension.   
\end{proof}
Let $k$ be a commutative ring below. Note that $\Psi$ is surjective iff the ideal $\{ e^1F(e^2) |\, F \in \hat{A},
e \in \mathcal{C} \} = R$. 
\begin{cor}
If $A \| B$ is a right free $k$-algebra extension of finite rank  where $A$ and $B$ are f.g.\ projective as 
$k$-modules,
and $\Psi$ is surjective, then $A \| B$ is a Frobenius extension.
\end{cor}
\begin{proof}
From the proof above we see that there is an epi $\vartheta:\, {}_BA_A \to {}_B \Hom (A_B, B_B)_A$
given by $a \mapsto Fa$ for some $F \in \hat{A}$.  But $\vartheta$ is an epi between
f.g.\ projective $k$-modules of the same $P$-rank with respect to localizations at prime
ideals $P$ in $k$, 
since $A_B$ is free.  
 It follows from a well-known fact that $\vartheta$ is bijective. (For similar reasons, 
a f.g.\ projective $k$-algebra is Dedekind-finite.)
Hence $A \| B$ is a free Frobenius extension. 
\end{proof}

\end{document}